\documentclass[number,seceqn,dvips]{arxbj}
\usepackage{cursive}  

\aid{0}
\volume{16}
\issue{4}
\pubyear{2010}
\firstpage{1224}
\lastpage{1239}
\doi{10.3150/10-BEJ249}

\makeatletter
\newremark{example}{Example}[section]
\newtheorem{theorem}{Theorem}[section]
\newproclaim{assumption}[theorem]{Assumption}
\newtheorem{lemma}[theorem]{Lemma}

\newremark{remark}[theorem]{Remark}
\newcommand{\dom}{\operatorname{dom}}
\newcommand{\interior}{\operatorname{int}}
\newcommand{\co}{\operatorname{co}}
\makeatother

\begin{document}
\begin{frontmatter}

\title{Testing composite hypotheses via convex~duality}
\runtitle{Testing composite hypotheses via convex duality}

\begin{aug}
\author[1]{\fnms{Birgit} \snm{Rudloff}\thanksref{1}\ead[label=e1]{brudloff@princeton.edu}\corref{}} \and
\author[2]{\fnms{Ioannis} \snm{Karatzas}\thanksref{2}\ead[label=e2]{ik@math.columbia.edu}}
\runauthor{B. Rudloff and I. Karatzas}
\address[1]{Department of Operations Research and Financial Engineering,
Princeton University, Princeton, NJ~08544, USA. \printead{e1}}
\address[2]{INTECH Investment Management, One Palmer
Square, Suite 441, Princeton, NJ 08542, USA and  Department of
Mathematics, Columbia University, MailCode 4438, New York, NY 10027,
USA.  \printead{e2}}
\end{aug}

\received{\smonth{6} \syear{2009}}
\revised{\smonth{11} \syear{2009}}

%
\begin{abstract}
We study the problem of testing composite hypotheses versus
composite alternatives, using a convex duality approach. In
contrast to classical results obtained by Krafft and Witting
(\textit{Z. Wahrsch. Verw.~Gebiete} \textbf{7} (1967) 289--302), where
sufficient optimality conditions are
derived via Lagrange duality, we obtain necessary and sufficient
optimality conditions via Fenchel duality under compactness
assumptions. This approach also differs from the methodology
developed in Cvitani\'c and Karatzas (\textit{Bernoulli} \textbf{7}
(2001) 79--97).
\end{abstract}

%

\end{frontmatter}
%

\section{Introduction}
The~problem of hypothesis testing is well understood in the
classical case of testing a simple hypothesis versus a simple
alternative. Suppose one wants to discriminate between two
probability measures $P$ (the ``null hypothesis'') and $Q$ (the
``alternative hypothesis''). In the classical Neyman--Pearson
formulation, one seeks a randomized test $ \varphi\dvtx  \Omega
\to[0,1] $ which is \textit{optimal}, in that it minimizes the
overall probability $  \mathbb{E}^Q [1 - \varphi] $ of not
rejecting $P$ when this hypothesis is false, while keeping below
a given significance level $\alpha\in(0,1)$ the overall
probability $  \mathbb{E}^P [ \varphi] $ of rejecting the
hypothesis $P$ when in fact it is true.

In this classical framework, an optimal randomized test
$ \widetilde{\varphi} $ always exists and can be calculated
explicitly in terms of a reference probability measure $R$, with
respect to which both measures are absolutely continuous (for
instance, $R= (P+Q)/2$). This test has the randomized 0--1
structure
%
\begin{equation}
\label{01}
\widetilde{\varphi} =  1_{\{ L >\mathfrak{z}\}} + \delta\cdot
1_{\{ L =\mathfrak{z}\}}
\end{equation}
which involves the likelihood ratio $ L = (\mathrm{d}Q/\mathrm{d}R)/(\mathrm{d}P/\mathrm{d}R) $ of
the densities of the null and the alternative hypotheses, the
quantile $  \mathfrak{z} = \inf\{ z \ge0 \dvt P ( L >z) \le
\alpha\} $ and the number $  \delta\in[0,1] $ which
enforces the significance-level requirement without slackness,
that is, $  \mathbb{E}^P [ \widetilde{\varphi}] = \alpha $.

The~problem becomes considerably more involved when the
hypotheses
are composite, that is, when one has to
discriminate between two entire \textit{families of probability
measures;} likelihood ratios of mixed strategies then have to be
considered. This type of problem also arises in the financial
mathematics context of minimizing the expected hedging loss in
incomplete or constrained markets; see, for example Cvitani\'c
\cite{Cvitanic00}, Schied \cite{Schied04} and Rudloff
\cite{Rudloff07}. It was shown by Lehmann \cite{Lehmann52}, Krafft
and Witting \cite{KafftWitting67}, Baumann \cite{Baumann68}, Huber
and Strassen \cite{HubStr73}, \"Osterreicher
\cite{Oesterreicher78}, Witting \cite{Witting85}, Vajda
\cite{Vajda89} and Cvitani\'c and Karatzas \cite{CviKar01}, that
duality plays a crucial role in solving the testing problem.
Most of these papers deal with Lagrange duality; they prove that
the typical 0--1 structure of (\ref{01}) is sufficient for
optimality and that it is both necessary and sufficient if a
dual solution exists. An important question then is to decide
when a dual solution will exist and to describe it when it does.

The~most recent of these papers, Cvitani\'c and Karatzas
\cite{CviKar01}, takes a different duality approach. Methods from
non-smooth convex analysis are employed and the set of densities
in the null hypothesis is enlarged in order to obtain the
existence of a dual solution, which again plays a crucial
role.

In the present paper, we shall use Fenchel duality. One advantage
of this approach is that as soon as one can prove the validity
of strong duality, the existence of a dual solution follows. We
shall show that strong duality holds under certain compactness
assumptions. This generalizes previous results insofar as no
need to enlarge the set of densities arises, a dual solution is
obtained and thus necessary and sufficient conditions for
optimality ensue.

In Section~\ref{sec2}, we introduce the problem of testing composite
hypotheses. Section~\ref{sec3} gives an overview of the duality results,
which are established and explained in detail in Section~\ref{sec4}. In
Section~\ref{sec5}, the imposed assumptions are discussed and possible
extensions are given. A comparison of the results and methods of
this paper with those in the existing literature can be found in
the last sections, notably Section~\ref{sec6}.

\section{Testing of composite hypotheses}\label{sec2}

Let $(\Omega,\mathcal{F})$ be a measurable space. A central
problem in the theory of hypothesis testing is to discriminate
between a given family $\mathcal P$ of probability measures
(composite ``null hypothesis'') and another given family $\mathcal
Q$ of probability
measures (composite ``alternative hypothesis'') on this space.

Suppose that there exists a \textit{reference
probability measure} $R$ on $(\Omega,\mathcal{F})$, that is, a
probability measure with respect to which all probability
measures $P \in\mathcal P$
and $Q\in\mathcal Q$ are absolutely continuous. We shall use
the notation $ Z_{\Pi} \equiv \mathrm{d} \Pi/\mathrm{d}R $ for the
Radon--Nikodym derivative of a finite measure $\Pi$ which is
absolutely continuous with respect to the reference measure and
$  \mathbb{E}^\Pi[Y]:= \int_\Omega Y \,\mathrm{d} \Pi= \int_\Omega Z_\Pi
Y \,\mathrm{d} R $ for the integral with respect to such $\Pi$ of an
$\mathcal{F}$-measurable function $Y\dvtx \Omega\rightarrow
[0,\infty)$. Finally, we shall denote the sets of these
Radon--Nikodym derivatives for the composite null hypothesis and
for the composite alternative hypothesis, respectively, by
\[
\mathfrak{Z}_{\mathcal P}:=\{Z_{P} \mid P\in\mathcal P\}
\quad \mbox{and}\quad
\mathfrak{Z}_{\mathcal Q}:=\{Z_Q \mid Q\in\mathcal Q\}.
\]
Both $\mathfrak{Z}_{\mathcal P} $ and $\mathfrak{Z}_{\mathcal Q}$
are subsets of the non-negative cone $  \mathbb{L}^1_+ $ and of
the unit ball in the Banach space $  \mathbb{L}^1\equiv
\mathbb{L}^1
(\Omega, \mathcal{F}, R)$. We assume
that
$   \Omega\times\mathfrak{Z}_\mathcal{P}
\ni(\omega, Z) \longmapsto
Z (\omega) \in[0,\infty) $ is
measurable with respect to the product $  \sigma$-algebra $
\mathcal{F} \otimes\mathcal{B} $, where $ \mathcal B $ is the
$\sigma$-algebra of Borel
subsets of $ \mathfrak{Z}_{\mathcal P}$.

We shall denote by $\Phi$ the set of all \textit{randomized tests,}
that is, of all Borel-measurable functions $\varphi\dvtx \Omega
\to[0,1]$ on $(\Omega,\mathcal{F})$. The~interpretation is as
follows: if the outcome $ \omega\in\Omega $ is observed and the
randomized test $  \varphi $ is used, then the null hypothesis
$ \mathcal P  $ is rejected with probability $ \varphi
(\omega)$. Thus, $ \mathbb{E}^P [\varphi]= \int_\Omega\varphi
(\omega) P (\mathrm{d} \omega)  $ is the overall probability of \textit{type I error} (of rejecting the null hypothesis, when in fact it
is true) under a scenario $  P \in\mathcal{P} $, whereas $
\mathbb{E}^Q[1-\varphi] $ is the overall probability of \textit{type II error} (of not rejecting the null hypothesis, when in
fact it is false) under the scenario $  Q \in\mathcal{Q} $.

We shall adopt the Neyman--Pearson point of view, whereby a
type I error is viewed as the more severe one and is not
allowed to occur with probability that exceeds a given acceptable
\textit{significance level} $\alpha\in(0,1)$, regardless of which
scenario $  P \in\mathcal{P} $ might materialize. Among all
randomized tests that observe this constraint,
%
\begin{equation}
\label{sign}
\mathfrak{s} (\varphi) :=  \sup_{P \in\mathcal P
} \mathbb{E}^{P }[\varphi] \leq \alpha,
\end{equation}
we then try to minimize the highest probability $  \sup_{Q \in
\mathcal Q } ( 1-\mathbb{E}^{Q }[\varphi]) $ of
type II error over all scenarios in the alternative hypothesis. In
other words, we look for a randomized test
$ \widetilde{\varphi} $ that maximizes the smallest \textit{power}
with respect to all alternative scenarios,
\[
  \pi(\varphi):= \inf_{Q\in\mathcal
Q}\mathbb{E}^Q[\varphi],
\]
over all randomized tests
$ \varphi $ whose `size' $\mathfrak{s} (\varphi) $, the quantity
defined in (\ref{sign}), does not exceed a given significance
level $\alpha$.

Equivalently, we look for a test $\widetilde{\varphi}
\in\Phi$ that attains the supremum
%
\begin{equation}\label{testop}
\underline{V}:=  \sup_{\varphi\in\Phi_\alpha} \pi(\varphi)
 =  \sup_{\varphi\in\Phi_\alpha} \Bigl( \inf_{Q\in\mathcal Q}
\mathbb{E}^Q[\varphi]\Bigr)
\end{equation}
of the power $ \pi(\varphi) $ over all generalized tests in
the
class
%
\begin{equation}\label{R_0}
\Phi_{\alpha}:=\Bigl\{\varphi\in\Phi \big|   \sup_{P \in
\mathcal{P }}
\mathbb{E}^{P}[\varphi]\leq
\alpha \Bigr\}.
\end{equation}
When such a randomized test $ \widetilde{\varphi} $ exists, it
will be called (max--min) \textit{optimal}.

\section{Duality}\label{sec3}

We shall denote by
$\Lambda_{+}$ the set of finite measures
on the measurable space $ (\mathfrak{Z}_{\mathcal P},\mathcal
B)$. We shall then associate to the maximization problem of
(\ref{testop}) the dual minimization problem
%
\begin{equation}
\label{opt pair test}
V^*:= \mathop{\inf_{ Q \in\mathcal Q }}_{\lambda\in\Lambda_{+}}
\mathcal{D} (Q, \lambda),
\end{equation}
where
%
\begin{equation}
\label{D}
\mathcal{D} (Q, \lambda):=
\mathbb{E}^R \biggl[\biggl(Z_Q-\int_{\mathfrak{Z}_{\mathcal P}}Z_{P }
\,\mathrm{d}\lambda
\biggr)^+ \biggr]
+\alpha  \lambda(\mathfrak{Z}_{\mathcal P}) .
\end{equation}
Here, and in the sequel, we view $  \int_{\mathfrak{Z}_{\mathcal
P}}Z_{P } (\omega) \,\mathrm{d}\lambda $ as the integral with respect to
the measure $  \lambda $ of the continuous functional $
\mathfrak{Z}_{\mathcal P} \ni Z \longmapsto\ell(Z; \omega):= Z
(\omega) \in\mathbb{R} $, for fixed $  \omega\in\Omega $; see
(\ref{pairing}) below for an amplification of this point.

The~idea behind (\ref{opt pair test}) and (\ref{D}) is
simple: we regard $  \lambda\in\Lambda_+ $ as a `Bayesian prior'
distribution
on the set $ \mathfrak{Z}_{\mathcal P} $ of densities for the null
hypothesis, whose effect is to reduce the composite null hypothesis $
\mathcal{P} $ to a simple one $  \{ P_*\} $ with $  Z_{ P_*} = \int
_{\mathfrak{Z}_{\mathcal P} }Z_{P }\,\mathrm{d} \lambda $ and whose total mass
$
\lambda(\mathfrak{Z}_{\mathcal P})<\infty $\vspace*{1pt} is a variable
that enforces the constraint in (\ref{sign}). More
precisely: for any given $  Q \in\mathcal{Q} $ and any $
\varphi\in\Phi_\alpha $, we have the \textit{weak duality}
\begin{eqnarray}
\label{dual}
\mathbb{E}^Q[\varphi] &=& \mathbb{E}^R[\varphi Z_Q] =
\mathbb{E}^R \biggl[\varphi\biggl(Z_Q-\int_{\mathfrak{Z}_{\mathcal P}
}Z_{P }\,\mathrm{d}\lambda\biggr)\biggr]
+\mathbb{E}^R \biggl[\varphi\int_{\mathfrak{Z}_{\mathcal P} }Z_{P
}\,\mathrm{d}\lambda\biggr]\nonumber
\\[-8pt]\\[-8pt]
&\le& \mathbb{E}^R
\biggl[\biggl(Z_Q-\int_{\mathfrak{Z}_{\mathcal P}}Z_{P } \,\mathrm{d}\lambda
\biggr)^+ \biggr]
+\alpha  \lambda(\mathfrak{Z}_{\mathcal P})  =  \mathcal{D} (Q,
\lambda), \qquad \lambda\in\Lambda_+ ,\nonumber
\end{eqnarray}
since
\begin{eqnarray*}
\mathbb{E}^R \biggl[\varphi\int_{\mathfrak{Z}_{\mathcal P} }Z_{P }\,\mathrm{d}
\lambda\biggr] & =&  \int_\Omega\varphi(\omega) \biggl( \int
_{\mathfrak{Z}_{\mathcal P}}Z_{P } (\omega)\, \mathrm{d}\lambda\biggr)\,   \mathrm{d}
R(\omega)\\
&  =&  \int_{\mathfrak{Z}_{\mathcal P}} \biggl( \int_\Omega
\varphi(\omega) Z_{P } (\omega) \, \mathrm{d} R(\omega) \biggr) \,  \mathrm{d}\lambda
\end{eqnarray*}
holds by the Fubini--Tonelli theorems, and this last quantity is dominated by
$  \alpha  \lambda( \mathfrak{Z}_{\mathcal P}) $ on the
strength of (\ref{sign}). We now observe that equality
holds in (\ref{dual}) if and only if both
%
\begin{eqnarray}\label{testB.a}
\varphi(\omega) =
 \cases{\displaystyle
1, &\quad  if $Z_{Q} (\omega)>\displaystyle\int_{\mathfrak{Z}_{\mathcal P}}Z_{P }(\omega)\, \mathrm{d}
\lambda$,
\vspace*{2pt}\cr
0, &\quad  if $ Z_{Q}(\omega)<\displaystyle\int_{\mathfrak{Z}_{\mathcal P}}Z_{P }
(\omega) \,\mathrm{d} \lambda$,}
\qquad   \mbox{for }
R\mbox{-a.e. } \omega\in\Omega
\end{eqnarray}
and
%
\begin{equation}\label{testB.b}
\mathbb{E}^{R} [ \varphi Z_P ]=\alpha,\qquad
 \mbox{for }
\lambda\mbox{-a.e. }  Z_P \in\mathfrak{Z}_{\mathcal
P}
\end{equation}
hold. It follows from (\ref{dual}) that the inequality $  \sup
_{\varphi
\in\Phi_\alpha} \mathbb{E}^Q[\varphi] \le\mathcal{D} (Q,
\lambda) $ holds for all $ \lambda\in\Lambda_+ $ and
$  Q \in\mathcal{Q} $ so that in the notation of (\ref{testop})
and (\ref{opt pair test}), we obtain
%
\begin{equation}
\label{dual.2}
\underline{V}    \le  \overline{V} :=  \inf_{Q\in\mathcal Q}
\Bigl(
\sup_{\varphi\in\Phi_\alpha} \mathbb{E}^Q[\varphi]\Bigr)  \le
  V^*.
\end{equation}

The~challenge, then, it to turn this `weak' duality into `strong' duality.
That is, to show that: equalities
$  \underline{V} = \overline{V}=V^* $ prevail in (\ref{dual.2});
the infimum in (\ref{opt pair test}) is attained by some\vspace*{1pt}
\mbox{$  (\widetilde{Q}, \widetilde{\lambda}) \in\mathcal Q \times
\Lambda_{+} $}; there exists a $  \widetilde{\varphi}\in
\Phi_\alpha $ for which the triple $  (\widetilde{\varphi},
\widetilde{Q},
\widetilde{\lambda}) $ satisfies (\ref{testB.a}), (\ref{testB.b});
for this triple equality prevails in (\ref{dual}); and the first element
$  \widetilde{\varphi} $ of this triple is optimal for the
generalized hypothesis-testing problem, that is, attains the supremum
in (\ref{testop}). We shall carry out this program, under appropriate
conditions, throughout the remainder of the paper.

When it exists, the measure $  \widetilde{\lambda} \in\Lambda_+ $
is called the ``least favorable distribution''. For explicit
computations of least favorable distributions in testing composite
hypotheses against simple alternatives (with $  \mathcal{Q} $ a
singleton), see Lehmann and Stein \cite{LehmannStein48}, Lehmann \cite
{Lehmann52} and Reinhardt \cite{Reinhardt61}, as well as Witting \cite
{Witting85}, pages 276--281 and Lehmann and Romano \cite
{LehmannRomano05}, Chapter 3. Here is an example abridged from these
last two sources.

\begin{example}
 Consider random variables $  X_1, X_2, \dots, X_n
 $ on $  (\Omega, \mathcal{F}) $ and probability measures on this
space under which these variables are independent with common Gaussian
distribution $  \mathcal{N} (\xi, \sigma^2) $. Thus, the random
variables $  \overline{X}= (1/n) \sum_{i=1}^n X_i $ and $ U=\sum
_{i=1}^n ( X_i - \overline{X} )^2 $ are sufficient for
the vector of parameters $(\xi, \sigma^2)$. For some given real
numbers $ \xi_1$ and $ \sigma_1^2 >0$, $\sigma_0^2 >0$, we shall
consider testing each of two composite hypotheses $\mathcal{P}^{(1)}$, $
\mathcal{P}^{(2)}$ versus the simple alternative $\mathcal{Q} = \{ Q\} $;
this latter corresponds to $( \xi, \sigma^2) = ( \xi_1, \sigma_1^2
)$. On the other hand, the hypothesis $  \mathcal{P}^{(1)} $ corresponds
to $  ( \xi, \sigma^2) \in\mathbb{R} \times[\sigma_0^2, \infty
) $, whereas the hypothesis $  \mathcal{P}^{(2)} $ corresponds to $  (
\xi, \sigma^2) \in\mathbb{R} \times(0, \sigma_0^2]$.

It is clear that the least favorable measure $  \widetilde{\lambda}
\in\Lambda_+$ should correspond to a distribution on $ \mathcal{B}
( \mathbb{R} \times(0, \infty))$ of the form $  \mu\otimes\delta
_{\sigma_0^2} $, where $  \mu $ is a measure on $ \mathcal{B} (
\mathbb{R})$; thus, under the measure $P_*$\vspace*{-1pt} with $  Z_{ P_*} = \int
_{\mathfrak{Z}_{\mathcal P} }Z_{P } \,\mathrm{d} \lambda $, the distribution
of $ \overline{X}$ has probability density
\[
\int_\mathbb{R} \bigl( 2 \curpi( \sigma_0^2 / n) \bigr)^{-1/2}  \exp
\biggl\{ - {   n (x-\xi)^2  \over2  \sigma_0^2 } \biggr\}  \mu(\mathrm{d}
\xi),  \qquad  x \in\mathbb{R},
\]
whereas, under the alternative $\mathcal{Q}$, the distribution of
$\overline{X}$ has probability density
\[
\bigl( 2 \curpi( \sigma_1^2 / n) \bigr)^{-1/2}  \exp\biggl\{ - {   n
(x-\xi_1)^2  \over2  \sigma_1^2 } \biggr\} ,\qquad    x \in\mathbb
{R}.
\]
\begin{itemize}
\item[$\bullet $] \textit{Testing $  \mathcal{P}^{(1)} $ versus $\mathcal{Q} $,
with} $  \sigma_0^2 > \sigma_1^2 $: The~least favorable $
\widetilde{\lambda} \in\Lambda_+$ corresponds to $  \mu= \delta
_{\xi_1} $. This way, the distribution of $\overline{X}$ is normal
under both $P_*$ and $Q$, with the same mean and with the smallest
possible difference in the variances.

\item[$\bullet $] \textit{Testing $  \mathcal{P}^{(2)} $ versus $\mathcal{Q} $,
with} $  \sigma_0^2 < \sigma_1^2 $: In this case, the least
favorable $  \widetilde{\lambda} \in\Lambda_+$ corresponds to $
\mu $ with Gaussian $  \mathcal{N} ( \xi_1, (\sigma_1^2 -
\sigma_0^2)/n) $ density. This guarantees that the distribution
of $\overline{X}$ is the same under both $P_*$ and $Q$.
\end{itemize}
\end{example}

\section{Results}\label{sec4}

In order to carry out the program outlined in the previous
section, we shall impose the following assumptions. A discussion
of their role can be found in Remark~\ref{Remark Assumpions}.

\begin{assumption}\label{Assumption}
\begin{enumerate}[(ii)]
\item[(i)] $ \mathfrak{Z}_{\mathcal Q} $ is a weakly compact,
convex subset of $ \mathbb{L}^1$.
\item[(ii)] $ \mathfrak{Z}_{\mathcal P}$ is a weakly compact subset
of $ \mathbb{L}^1$.
\end{enumerate}
\end{assumption}

Our main result reads as follows.

\begin{theorem}[(Generalized Neyman--Pearson lemma)]\label{Theorem Testtheorie}
Let $ \mathcal P, \mathcal Q $ be families of probability
measures on $  (\Omega, \mathcal{F})$, as in Sections
\textup{\ref{sec2}} and \textup{\ref{sec3}}, that satisfy
Assumption~\textup{\ref{Assumption}}.
For a given constant $\alpha\in(0,1)$, recall the subclass
$  \Phi_\alpha $ of randomized tests in \textup{(\ref{R_0})}.

There then exists a randomized test $ \widetilde{\varphi} \in
\Phi_\alpha $
which
attains the supremum in
\textup{(\ref{testop})}. There also exists a solution to the dual problem
of \textup{(\ref{opt pair test})}, namely, a pair
$(\widetilde{Q},\widetilde{\lambda})\in\mathcal Q\times
\Lambda_{+}$ which attains the infimum there.

Furthermore, strong duality is satisfied, in the sense that:
\begin{itemize}
\item the optimal test for (\ref{testop}) has
the structure of (\ref{testB.a}), (\ref{testB.b}),
namely
%
\begin{eqnarray}
\widetilde{\varphi} (\omega) = \cases{
1, &\quad if $Z_{\widetilde{Q}} (\omega)>\displaystyle\int_{\mathfrak{Z}_{\mathcal P}}Z_{P
}(\omega) \,\mathrm{d}\widetilde{\lambda}$,\vspace*{2pt}
\cr
0, &\quad if $ Z_{\widetilde{Q}}(\omega)<\displaystyle\int_{\mathfrak{Z}_{\mathcal P}}Z_{P
}(\omega) \,\mathrm{d}\widetilde{\lambda}$,}
\qquad   \mbox{for }
R\mbox{-a.e. } \omega\in\Omega
\end{eqnarray}
and
%
\begin{equation}\label{testB}
\mathbb{E}^{R}[\widetilde{\varphi} Z_P ]=\alpha \qquad
 \mbox{for }
\widetilde{\lambda}\mbox{-a.e. } Z_P \in\mathfrak{Z}_{\mathcal P},
\mbox{ whereas}
\end{equation}
\item$(\widetilde{\varphi},\widetilde{Q})$ is a
saddlepoint in $
\Phi_\alpha\times
\mathcal{Q} $ of the functional
$ (\varphi,Q)\mapsto\mathbb{E}^{Q}[\varphi] $, namely,
%
\begin{equation}
\label{dual.3}
\mathbb{E}^{\widetilde{Q}}[\varphi] \le
\mathbb{E}^{\widetilde{Q}}[\widetilde{\varphi} ] \le
\mathbb{E}^Q[\widetilde{\varphi} ]\qquad    \forall ( \varphi,
Q) \in  \Phi_\alpha\times\mathcal{Q}.
\end{equation}
\end{itemize}
\end{theorem}

The~theorem will be proven in several steps, using the
following lemmata. The~first of these, Lemma \ref{lemma}, seems to be
well known (cf.~the Appendix in Lehmann and Romano \cite
{LehmannRomano05}). We could not, however, find in the literature a
result in exactly the form we needed that we could cite directly, so we
provide a proof for the sake of completeness. We shall freely use the
convention of denoting by ``max'' (resp., ``min'') a supremum (resp.,
infimum) which is attained.

\begin{lemma}\label{lemma}
The~supremum in \textup{(\ref{testop})} is attained by some
randomized test $ \widetilde{\varphi} \in\Phi_\alpha $ and
there exists a $  \widetilde{Q} \in\mathcal{Q} $ such that the
saddlepoint property
\textup{(\ref{dual.3})} holds. In particular, the lower- and upper-values $
\underline{V} $ and
$ \overline{V} $ of
\textup{(\ref{testop})} and \textup{(\ref{dual.2})}, respectively, are the same:
%
\begin{equation}\label{saddle}
\max_{\varphi\in\Phi_{\alpha}} \Bigl( \min_{Q\in\mathcal
Q}\mathbb{E}^Q[\varphi]\Bigr)
 = \min_{Q\in\mathcal Q} \Bigl( \max_{\varphi\in
\Phi_{\alpha}}\mathbb{E}^Q[\varphi]\Bigr).
\end{equation}
\end{lemma}

\begin{pf}
It is well known that the set $  \Phi $ of all randomized tests
is weakly-* compact (this follows from weak sequential compactness;
cf.~\cite{Witting85}, page 270 or \cite{LehmannRomano05},
Theorem~A.5.1). We give a short proof for the sake of completeness. The~subset $  \Phi $ of the Banach space $
\mathbb{L}^\infty\equiv\mathbb{L}^\infty(\Omega, \mathcal{F},
R) $ is weakly-* compact, as it is a weakly-* closed subset of the
weakly-* compact
unit ball in $\mathbb{L}^\infty$ (Alaoglu's theorem; see, e.g., Dunford
and Schwartz \cite{DunSchw88}, Theorem~V.4.2 and Corollary V.4.3).
To see that $  \Phi $ is weakly-* closed, consider a net
$\{\varphi_\alpha\}_{\alpha\in D}\subseteq\Phi$ that converges to
$\varphi$ with respect to the weak-* topology in
$\mathbb{L}^\infty$. This means that for all $X\in\mathbb{L}^1$,
we have $\mathbb{E}^R[\varphi_\alpha X]\rightarrow
\mathbb{E}^R[\varphi X]$. If there existed an event
$\Omega_1\in\mathcal F$ with $R(\Omega_1)>0$ and $  \{
\varphi>1\} \subseteq\Omega_1 $, then we could choose
$\widehat{X}(\omega)=1_{\Omega_1}(\omega)\in\mathbb{L}^1$ and
obtain $\mathbb{E}^R[\varphi\widehat{X}]>R(\Omega_1)$. However, this
contradicts
$\mathbb{E}^R[\varphi\widehat{X}]=\lim_{\alpha}\mathbb
{E}^R[\varphi_\alpha
\widehat{X}]\leq R(\Omega_1)$, which follows from
$\varphi_\alpha\leq1$ for all $\alpha\in D$ since
$\varphi_\alpha\in\Phi$. Hence, $\varphi\leq1 $ holds $R$-a.e.
It can be similarly shown that $\varphi\geq0$ also holds
$  R$-a.e.

Thus, $\Phi$ is indeed weakly-* closed, hence weakly-* compact.
Since the mapping $ \varphi\mapsto\sup_{P \in\mathcal{P}
}\mathbb{E}^{P
}[\varphi] $ is lower-semicontinuous in the weak-* topology, the
set $\Phi_\alpha$ in (\ref{R_0}) is \mbox{weakly-*} closed, hence weakly-*
compact. Because of the upper-semicontinuity of the mapping
$ \varphi\mapsto\pi(\varphi) =
\inf_{Q\in\mathcal{Q}}\mathbb{E}^{Q}[\varphi] $ in the weak-*
topology, there exists a $ \widetilde{\varphi}\in\Phi_\alpha $
that attains the supremum in (\ref{testop}).

The~weak-* compactness and convexity of $ \Phi_\alpha $ together
with the
weak compactness and convexity of $\mathfrak{Z}_{\mathcal Q}$
(Assumption~\ref{Assumption}(i)) enable us to apply the von
Neumann/Sion minimax theorem (see, for instance, \cite{Aubin79},
Theorem~7, Section~7.1, or \cite{Aubin00}, Section 2.7,
pages 39--45); the assertions follow.
\end{pf}

Let us now fix an arbitrary $Q\in\mathcal Q $ and consider as
our primal problem the inner maximization in the middle term of
(\ref{dual.2}), namely,
%
\begin{equation}\label{ipNP}
\mathfrak{p}(Q):=
\sup_{\varphi\in\Phi_\alpha}\mathbb{E}^Q[\varphi].
\end{equation}
This supremum is always attained since
$ \Phi_\alpha $ is weakly-* compact. We want to show that strong
duality holds
between (\ref{ipNP}) and its \textit{Fenchel dual problem} which,
we claim, is of the form
%
\begin{equation}
\label{lambda}
  \mathfrak{d}(Q) = \inf_{\lambda\in\Lambda_{+}} \mathcal{D} (Q,
\lambda) = \inf_{\lambda\in\Lambda_{+}}\biggl[ \int
_{\Omega}
\biggl(Z_Q-\int _{\mathfrak{Z}_{\mathcal P} }Z_{P
} \,\mathrm{d}\lambda\biggr)^{+}\,\mathrm{d}R
 + \alpha \lambda(\mathfrak{Z}_{\mathcal P})\biggr].
\end{equation}
Note that, from (\ref{dual}), we have $  \mathfrak{p}(Q) \le
\mathfrak{d}(Q)$.

In this setting, the typical 0--1 structure of the randomized
test $ \widetilde{\varphi}_Q \in\Phi_\alpha $ that attains the
supremum in (\ref{ipNP}) is necessary and sufficient for
optimality.

\begin{lemma}\label{theorem opt Test NP}
Strong duality holds for problems \textup{(\ref{ipNP})} and
\textup{(\ref{lambda})}, that is,
\[
\forall Q\in\mathcal Q,\qquad
\mathfrak{d}(Q) = \mathfrak{p}(Q).
\]
Moreover, for each $ Q\in\mathcal Q $, there exists a measure
$ \widetilde{\lambda}_{Q} \in\Lambda_+ $
which attains
the infimum in \textup{(\ref{lambda})},
whereas an optimal test
$ \widetilde{\varphi}_Q \in\Phi_\alpha $ that attains the
supremum in \textup{(\ref{ipNP})} exists and has
the structure of \textup{(\ref{testB.a})} and \textup{(\ref{testB.b})}, namely,
%
\begin{eqnarray}
\label{NP Loesung 1}
\widetilde{\varphi}_Q (\omega) = \cases{
1, &\quad if $Z_{Q} (\omega)>\displaystyle\int_{\mathfrak{Z}_{\mathcal P}}Z_{P }
(\omega) \,\mathrm{d}\widetilde{\lambda}_Q$,\vspace*{2pt}
\cr
0, &\quad if $Z_{Q}(\omega)<\displaystyle\int_{\mathfrak{Z}_{\mathcal P}}Z_{P }
(\omega) \,\mathrm{d}\widetilde{\lambda}_Q$,}
\qquad  \mbox{for }
R\mbox{-a.e. } \omega\in\Omega
\end{eqnarray}
and
%
\begin{equation}
\label{NP Loesung 2}
\mathbb{E}^{R} [ \widetilde{\varphi}_Q Z_P ] = \alpha
\qquad  \mbox{for }
\widetilde{\lambda}_Q\mbox{-a.e. } Z_P \in\mathfrak{Z}_{\mathcal
P}.
\end{equation}
\end{lemma}

\begin{pf}
Let $ \mathcal L $ be the linear space of all continuous
functionals $ \ell\dvtx \mathfrak{Z}_{\mathcal{P}}\rightarrow\mathbb
{R} $ on
the weakly compact subset $ \mathfrak{Z}_{\mathcal{P} }$ of $
\mathbb{L}^1 $ (Assumption~\ref{Assumption}(ii)) with pointwise
addition and multiplication by real numbers and pointwise partial
order
%
\begin{equation}
\label{ordnung}
\ell_1\leq\ell_2 \quad  \Longleftrightarrow\quad    \ell_2-\ell_1  \in
  \mathcal L_+:= \{\ell\in\mathcal
L  \mid   \ell(Z_{P})\geq0, \forall P \in
\mathcal{P}\}.
\end{equation}
We endow $ \mathcal L $ with the supremum
norm $ \|\ell\|_\mathcal L=\sup_{P\in
\mathcal{P}}|\ell(Z_{P})| $, which ensures that $ \mathcal L $
is a Banach space (Dunford and Schwartz \cite{DunSchw88},
Section~IV.6).

Similarly, we let $ \Lambda $ be the space of regular finite signed
measures $  \lambda= \lambda^+ - \lambda^- $ on
$ (\mathfrak{Z}_{\mathcal{P} }, \mathcal B) $, with $
\lambda^\pm\in\Lambda_+ $. The~space $ \Lambda $ is sometimes
denoted $\mathfrak{ca}_r(\mathfrak{Z}_{\mathcal{P} }, \mathcal B
)$ in the literature, in order to stress that it consists of \textit{r}egular $\mathfrak{ c}$ountably $\mathfrak{a}$dditive signed
measures of finite variation (\cite{AliBor06}, Definition~12.11). We
regard this space as the
norm-dual of $ \mathcal L $ with the bilinear form
%
\begin{equation}
\label{pairing}
\langle\ell, \lambda\rangle = \int_{\mathfrak{Z}_{\mathcal P}
}\ell \,\mathrm{d}\lambda  \qquad  \mbox{for }   \ell\in\mathcal L,
\lambda\in\Lambda;
\end{equation}
see Aliprantis and Border \cite{AliBor06}, Corollary~14.15.
(Intuitively speaking, the elements of $ \Lambda $ are
generalized Bayesian priors that may assign negative mass to
certain null hypotheses.)
We also note the clear identification $  \Lambda_{+} \equiv
\{\lambda\in\Lambda \mid   \lambda(B)\geq0,  \forall
B\in\mathcal B\} $.

  Let us define a linear
operator
$ A\dvtx (\mathbb{L}^\infty,\|\cdot\|_{\mathbb{L}^\infty})
\to(\mathcal L,\|\cdot\|_\mathcal L) $ by
%
\begin{equation}
\label{OperA}
\mathbb{L}^\infty \ni\varphi \quad  \longmapsto \quad  (A\varphi)(Z_{P})
:= -\mathbb{E}^{P}[\varphi] =  - \mathbb{E}^R [
\varphi Z_P]  \in\mathbb{R}
\end{equation}
for $Z_P \in\mathfrak{Z}_\mathcal{P}$. This operator is bounded,
thus continuous. We also introduce the constant functionals
$ \mathbf{1},\mathbf{0} \in\mathcal L$ by
\[
\forall Z_P \in
\mathfrak{Z}_\mathcal{P},\qquad \mathbf{1}(Z_{P})=1\in\mathbb{R}, \qquad
\mathbf{0}(Z_{P})=0\in\mathbb{R}.
\]
The~constraint of (\ref{sign}) can then be rewritten as
\[
\alpha\mathbf{1}+A\varphi\geq\mathbf{0}
\quad  \Longleftrightarrow\quad
A\varphi\in\mathcal{L}_{+}-\alpha\mathbf{1}.
\]
With this notation, for any given $  Q \in\mathcal{Q},$ the
primal problem (\ref{ipNP}) is cast as
\begin{eqnarray}
\label{pi neu}
-\mathfrak{p}(-Q)
 &=& \inf_{\varphi\in\mathbb{L}^\infty}
\bigl(  \mathbb{E}^Q[\varphi]
+\mathcal I_{\Phi}(\varphi)
+\mathcal I_{
\mathcal{L}_+-\alpha\mathbf{1}}
(A\varphi)\bigr)\nonumber
\\[-8pt]\\[-8pt]
\nonumber &=& \inf_{\varphi\in\mathbb{L}^\infty}
\bigl(  f (\varphi)+ g (A\varphi)\bigr)
,
\end{eqnarray}
where $  -Q $ is interpreted as a finite, signed measure on $
(\Omega, \mathcal{F}) $ (cf.~the discussion preceding
(\ref{pairing})).
Here, and for the remainder of this proof, we use the notation $
\mathcal{I}_C (\varphi):= 0 $ for $  \varphi\in C$, $
\mathcal{I}_C (\varphi):= \infty$ for $  \varphi\notin C$,
as
well as
%
\begin{equation}
\label{FG}
f(\varphi):=\mathbb{E}^Q[\varphi]+\mathcal I_{\Phi}(\varphi),\qquad
  g(A\varphi):=\mathcal I_{\mathcal
L_+-\alpha\mathbf{1}}(A\varphi).
\end{equation}

 \noindent$\bullet$ We claim that the \textit{Fenchel dual} of the primal problem in (\ref{ipNP}) has the form
(\ref{lambda}). We shall begin the proof of this claim by
recalling from Ekeland and Temam \cite{EkeTem76},
Proposition~III.1.1, Theorem~III.4.1 and Remark~III.4.2 that the
Fenchel dual of the problem (\ref{pi neu}) is given by
%
\begin{equation}
\label{d inneres}
-\mathfrak{d}(-Q) = \sup_{\lambda\in\Lambda}
\bigl(-f^{*}(A^{*}\lambda)-
g^{*}(-\lambda)\bigr),
\end{equation}
where $ A^{*}\dvtx  \Lambda\to\mathfrak{ba}(\Omega,\mathcal{F},R) $
is the adjoint of the operator $A$ in (\ref{OperA}). Here, and in
the sequel, $  \mathfrak{ba}(\Omega,\mathcal{F}, R) $ is the
space of $\mathfrak{b}$ounded, (finitely-) $\mathfrak{a}$dditive
set-functions on $
(\Omega,\mathcal{F}) $ which are absolutely continuous with respect
to $ R $; see, for instance, Yosida \cite{Yosida95}, Chapter IV,
Section~9, Example 5.

The~function $g^*(\cdot)$ is the conjugate of the function
$g(\cdot)$, namely,
\begin{eqnarray*}
g^{*}(\lambda)&=&\sup_{\widetilde{\ell}\in\mathcal L}
\bigl(\langle \widetilde{\ell},\lambda\rangle-
\mathcal I_{\mathcal L_{+}
-\alpha\mathbf{1}}(\widetilde{\ell} )\bigr)
 = \sup_{\widetilde{\ell} \in \mathcal L_{+}-\alpha\mathbf{1}}
\langle
\widetilde{\ell},\lambda\rangle
 = \sup_{\ell\in\mathcal L_{+}}\langle\ell-\alpha\mathbf{1},
\lambda\rangle
\\
&=&\sup_{\ell\in\mathcal L_{+}}\langle
\ell,\lambda\rangle-\alpha\int _{\mathfrak{Z}_{\mathcal P} }
\mathrm{d}\lambda
 = \mathcal I_{\mathcal
L_{+}^{*}}(\lambda)-\alpha \lambda(\mathfrak{Z}_{\mathcal{P} }),
\end{eqnarray*}
where
%
\begin{equation}
\label{StarPlus}
\mathcal L_{+}^{*}:= \{\lambda\in\Lambda \mid
 \langle\ell,\lambda \rangle\leq0,  \forall
\ell\in\mathcal
L_+ \}
\end{equation}
is the negative dual cone of $ \mathcal L_{+} $. The~last equality in the above string holds because the set $ \mathcal
L_{+} $ defined in (\ref{ordnung}) is a cone containing the origin $
\mathbf{0}\in\mathcal L $.

To determine the conjugate $f^*(\cdot)$ of the function $f(\cdot)$
at $A^{*}\lambda$, namely
\[
f^{*}(A^{*}\lambda) = \sup_{\varphi\in\mathbb{L}^\infty}\{
\ensuremath{ \langle A^{*}\lambda,\varphi \rangle}
-\mathbb{E}^Q[\varphi]-\mathcal I_{\Phi}(\varphi)\},
\]
we have to calculate $\langle A^{*}\lambda,\varphi\rangle$.
By the definition of $A^*$, the equation $\langle
A^{*}\lambda,\varphi\rangle=\langle\lambda,A\varphi\rangle$ has to
be satisfied for all $ \varphi\in
\mathbb{L}^{\infty}, \lambda\in\Lambda $ (see \cite{AliBor06},
Chapter~6.8). Thus,
\[
\forall  \varphi\in\mathbb{L}^{\infty}, \forall  \lambda
\in\Lambda,\qquad
 \langle A^{*}\lambda,\varphi\rangle =  -
\int _{\mathfrak{Z}_{\mathcal P}}\mathbb{E}^{R}
[\varphi Z_P] \,\mathrm{d}\lambda
\]
and the conjugate of the function $f(\cdot)$ at $A^{*}\lambda$
is evaluated as
\[
f^{*}(A^{*}\lambda)
 = \sup_{\varphi\in\Phi} \biggl( -\int _{ \mathfrak
{Z}_{\mathcal P} }
\mathbb{E}^{R}[\varphi Z_P]\,\mathrm{d}\lambda-\mathbb{E}^Q[\varphi]  \biggr).
\]
The~dual problem (\ref{d inneres}) therefore
becomes
\begin{eqnarray}\label{dproof2-}
-\mathfrak{d}(-Q)&=&\sup_{\lambda\in\Lambda}
\biggl[ -\sup_{\varphi\in\Phi}\biggl(-\int _{ \mathfrak
{Z}_{\mathcal P} }\mathbb{E}^{R }
[\varphi Z_P]\,\mathrm{d}\lambda-\mathbb{E}^Q[\varphi]\biggr)
-\mathcal I_{ -\mathcal L_{+}^{*}}(\lambda)
-\alpha \lambda(\mathfrak{Z}_{\mathcal{P}})\biggr]\nonumber
\\[-8pt]\\[-8pt]
&=&\sup_{\lambda\in-\mathcal L_{+}^{*}}
\biggl[ -\sup_{\varphi\in\Phi}\biggl(-
\int _{\mathfrak{Z}_{\mathcal P}}\mathbb{E}^{R}
[\varphi Z_P]\,\mathrm{d}\lambda
-\mathbb{E}^Q[\varphi]\biggr)-\alpha
\lambda(\mathfrak{Z}_{\mathcal{P}}) \biggr].\nonumber
\end{eqnarray}
It is not hard to show the property $ -\mathcal L_{+}^{*}=
\Lambda_{+} $ for the set in (\ref{StarPlus}),
so the expression of (\ref{dproof2-}) can be recast in the form
%
\begin{equation}
\label{dproof2}
-\mathfrak{d}(-Q) = \sup_{\lambda\in\Lambda_{+}}
\biggl[ -\sup_{\varphi\in\Phi}\biggl(-
\int _{\mathfrak{Z}_{\mathcal P} }
\mathbb{E}^R[Z_{P
} \varphi] \,\mathrm{d}\lambda-\mathbb{E}^R[Z_Q \varphi]\biggr)
-\alpha \lambda(\mathfrak{Z}_{\mathcal{P} }) \biggr].
\end{equation}

\noindent$\bullet$ Now, both $(\Omega,\mathcal{F},R)$ and
$(\mathfrak{Z}_{\mathcal{P} }, \mathcal B,\lambda)$ for
$\lambda\in\Lambda_{+}$ are positive, finite measure spaces.
The~mapping $   \Omega\times
\mathfrak{Z}_\mathcal{P}  \ni(\omega, Z_P) \longmapsto
f(\omega,Z_P):=Z_P (\omega)\varphi(\omega) \in\mathbb{R} $ is
measurable with respect to the product $  \sigma
$-algebra $
\mathcal{F} \otimes\mathcal{B} $ for every
$ \varphi\in\Phi$, thanks to the measurability
assumption of Section~\ref{sec2}, whereas
\[
\int _{\mathfrak{Z}_{\mathcal P}}\int _\Omega|Z_{P }
\varphi|  \,\mathrm{d}R\,
\mathrm{d}\lambda  \leq  \Bigl(\sup_{P \in\mathcal{P} }\Vert Z_{P
}\Vert _{ \mathbb{L}^1}
\Bigr)\lambda(\mathfrak{Z}_{P })
 = \lambda(\mathfrak{Z}_{\mathcal{P} }) <  \infty
\]
holds for every $ \lambda\in\Lambda_{+} $ and $ \varphi\in
\Phi $ since $ \|\varphi\|_{  \mathbb{L}_\infty}\leq1 $. Thus,
we can
apply the Fubini--Tonelli theorem (see \cite{DunSchw88},
Corollary~III.11.15) and deduce that the order of integration can
be interchanged, that is, for all $\lambda\in\Lambda_{+}$ and all
$\varphi\in\Phi$, we have
\[
\label{tonelli}
\int _{\mathfrak{Z}_{\mathcal P} }\int _\Omega Z_{P
}\varphi \,\mathrm{d}R \,\mathrm{d}\lambda
 = \int _\Omega\int _{\mathfrak{Z}_{\mathcal
P}}Z_{P }\varphi \, \mathrm{d}\lambda \,\mathrm{d}R  < \infty.
\]

In (\ref{dproof2}), only elements $\lambda\in\Lambda_{+}$ and
$\varphi\in\Phi$ are considered, so we can interchange the order
of integration and obtain
%
\begin{equation}
\label{dproof2b}
-\mathfrak{d}(-Q) = \sup_{\lambda\in\Lambda_{+}}\biggl(-\sup
_{\varphi\in
\Phi}
\mathbb{E}^R \biggl[ \varphi\biggl(-Z_Q-\int _{\mathfrak
{Z}_{\mathcal P} }
Z_{P } \,\mathrm{d}\lambda\biggr)  \biggr]-\alpha\lambda(\mathfrak
{Z}_{\mathcal{P}
})\biggr).
\end{equation}
Since $ \varphi\in\Phi $ is a randomized test,
the supremum over all $ \varphi\in\Phi $ in (\ref{dproof2b}) is
attained by some $ \overline{\varphi}_{\lambda, -Q} \in\Phi $
of a form similar to (\ref{testB.a}), namely,
%
\begin{equation}
\label{varphi quer}
\overline{\varphi}_{\lambda, -Q}
(\omega) = \cases{
1, &\quad if $ - Z_{Q} (\omega)>
\displaystyle\int_{\mathfrak{Z}_{\mathcal P}}Z_{P }(\omega) \,\mathrm{d} \lambda$\vspace*{2pt}
\cr
0, &\quad  if $ - Z_{Q}(\omega)<
\displaystyle\int_{\mathfrak{Z}_{\mathcal P}}Z_{P }(\omega) \,\mathrm{d} \lambda$}
\qquad  \mbox{for }
R\mbox{-a.e. } \omega\in\Omega.
\end{equation}

Given any finite, signed measure $  \Pi= \Pi^+ - \Pi^- $ on $
(\Omega, \mathcal{F}) $ with $  \Pi^\pm\ll R $, let us write
$  Z_\Pi= Z_{\Pi^+} - Z_{\Pi^-}$, define
%
\begin{equation}
\label{Ypsilon}
\Upsilon_{\lambda, \Pi} :=  Z_\Pi-\int_{\mathfrak{Z}_{\mathcal
P} }Z_{P }\,\mathrm{d}\lambda  \in \mathbb{L}^1
\end{equation}
and let
$ \Upsilon_{\lambda, \Pi}^{+}$ (resp.,
$ \Upsilon_{\lambda, \Pi}^{-} $) be the positive
(resp., negative) part of the random variable in
(\ref{Ypsilon}).
With this notation, and recalling (\ref{varphi quer}), the value
of the dual problem (\ref{dproof2b}) becomes
\[
-\mathfrak{d}(-Q) = \sup_{\lambda\in\Lambda_{+}}
\{-\mathbb{E}^R
[\Upsilon_{\lambda, -Q}^{+} ]-\alpha\lambda( {\mathfrak
{Z}_{\mathcal P}})
\}
\]
and thus
%
\begin{equation}
\label{pLUSmINUS}
\mathfrak{d}( Q) =
\inf_{\lambda\in\Lambda_{+}}   \{ \mathbb{E}^R
[\Upsilon_{\lambda, Q}^{+} ]+\alpha\lambda( {\mathfrak
{Z}_{\mathcal P}})
\}.
\end{equation}
We deduce from this representation and (\ref{Ypsilon}) that the
dual $\mathfrak{d}(Q)$ of the primal problem $\mathfrak{p}(Q)$ of
(\ref{ipNP}) is indeed as claimed in equation (\ref{lambda}).

 \noindent$\bullet$
Strong duality now holds if both $f(\cdot)$ and $g(\cdot)$ are
convex, $g(\cdot)$ is continuous at some $A\varphi_{0}$ with
$\varphi_{0}\in\dom(f)$ and $\mathfrak{p}(Q)$ is finite (see
\cite{EkeTem76}, Theorem~III.4.1 and
Remark~III.4.2).

Indeed, $\mathfrak{p}(Q)$ is clearly finite. In (\ref{FG}), the function
$f(\cdot)$ is
convex since $\Phi$ is a convex set; the function $g(\cdot)$ is convex
since the set $ \mathcal L_{+}-\alpha\mathbf{1} $ is convex, as
well as
continuous at some $A\varphi_{0}$
with $\varphi_{0}\in\dom( f)$, provided that
$A\varphi_{0}\in\interior(\mathcal L_{+}-\alpha\mathbf{1})$. If we
take $\varphi_{0}\equiv0$, then $ \varphi_0\in\dom(f) $ since
$\varphi_0\in\Phi$, and we see that
$A\varphi_{0}=\mathbf{0}\in\interior(\mathcal
L_{+}-\alpha\mathbf{1})$ since $\interior( \mathcal L_{+})
\neq\varnothing $ in the norm topology and $\alpha>0$. Hence, we
have strong duality.

\noindent$\bullet$
The~existence of a
solution to the primal problem $\mathfrak{p}(Q)$ (i.e., of a
generalized test $ \widetilde{\varphi}_Q\in\Phi_\alpha $ which
attains the supremum in (\ref{ipNP})) follows from the weak-*
compactness of $\Phi_\alpha$. With strong duality established, the
existence of a solution to the dual problem, that is, of an element
$ \widetilde{\lambda}_Q\in\Lambda_{+} $ which attains the infimum
in (\ref{lambda}), follows (\cite{EkeTem76},
Theorem~III.4.1 and
Remark~III.4.2), whereas the values of the primal (resp., the
dual) objective functions at $\widetilde{\varphi}_Q$ (resp.,
$\widetilde{\lambda}_Q$) coincide. To indicate the dependence of
these quantities on the selected $ Q\in\mathcal Q $, we have used
the notation $\widetilde{\varphi}_Q$ and $\widetilde{\lambda}_Q$ for
the primal and dual solutions, respectively.

 \noindent$\bullet$
\textit{These considerations lead to a necessary and sufficient
condition for optimality}. Indeed, let us write the expression for $
\mathbb{E}^Q[\varphi] $ from (\ref{dual}) as
\[
\mathbb{E}^Q[\varphi] = \mathbb{E}^R[\varphi\Upsilon_{\lambda,
Q}^{+} ]-
\mathbb{E}^R[\varphi\Upsilon_{\lambda, Q}^{-} ]
+\mathbb{E}^R \biggl[\varphi\int_{\mathfrak{Z}_{\mathcal P} }
Z_{P }\,\mathrm{d}\lambda\biggr],
\]
in the notation of (\ref{Ypsilon}), and subtract it from
the dual objective function $  \mathbb{E}^R[ \Upsilon_{\lambda,
Q}^{+} ] + \alpha\lambda(\mathfrak{Z}_\mathcal{P})$, as in
(\ref{pLUSmINUS}). Because of strong duality, this difference has
to be zero when evaluated at $  (\varphi, \lambda) =
(\widetilde{\varphi}_Q,
\widetilde{\lambda}_Q) $, namely,
\[
\mathbb{E}^R [\Upsilon_{\widetilde{\lambda}_Q, Q}^{+}
(1-\widetilde{\varphi}_Q)]
+\mathbb{E}^R [ \Upsilon_{\widetilde{\lambda}_Q, Q}^{-}
 \widetilde{\varphi}_Q ]
+\int _{\mathfrak{Z}_{\mathcal P}}(\alpha-\mathbb{E}^R
[Z_{P } \widetilde{\varphi}_Q])\, \mathrm{d}\widetilde{\lambda}_Q =
0.
\]
Each of these three integrals is non-negative, so their sum is
zero if and only if $ \widetilde{\varphi}_Q\in\Phi_\alpha $
satisfies the condition
(\ref{NP Loesung 2}) of Lemma~\ref{theorem opt Test NP} and is of
the form (\ref{NP Loesung 1}) or, equivalently, of the form $
\widetilde{\varphi}_Q \equiv
\overline{\varphi}_{ \widetilde{\lambda}_Q, -Q} $ of
(\ref{varphi quer}).
\end{pf}

We are ready now to prove our main result.

\begin{pf*}{Proof of Theorem~\ref{Theorem Testtheorie}}
Lemma~\ref{theorem opt Test NP} guarantees that
\[
\min _{Q\in\mathcal{Q}}\Bigl( \max _{\varphi\in
\Phi_\alpha} \mathbb{E}^{Q}[\varphi]\Bigr)
  =  \min _{Q\in\mathcal{Q}} \mathfrak{p}(Q) =
\min _{Q\in
\mathcal{Q}} \mathfrak{d}(Q) =   \min_{(Q, \lambda) \in
\mathcal\mathcal{Q} \times\Lambda_{+}}
\mathcal{D} (Q, \lambda)  =  V^*
\]
in the notation of (\ref{opt pair test}) and (\ref{ipNP}),
(\ref{lambda}). From Lemma~\ref{lemma}, it follows that there
exists an element $ \widetilde{Q}  $ of $ \mathcal{Q} $ which
attains the infimum in (\ref{dual.2}).
For this $ \widetilde{Q}  $, Lemma~\ref{theorem opt Test NP}
shows the existence of an element
$ \widetilde{\lambda}_{\widetilde{Q}} $ of $  \Lambda_+ $ that attains
the infimum in (\ref{lambda}). Thus, there
exists a pair $(\widetilde{Q},\widetilde{\lambda})$ that attains
the infimum in (\ref{opt pair test}) and Lemma~\ref{theorem opt
Test NP} gives the required structural result.
\end{pf*}

\begin{corollary*}
It follows that the optimal
randomized test has the form
%
\begin{equation}
\label{012}
\widetilde{\varphi}(\omega) = 1_{\{Z_{\widetilde{Q}} >
\int_{\mathfrak{Z}_{\mathcal P} }Z_{P } \mathrm{d} \widetilde{\lambda}
\}} (\omega)
 + \delta(\omega) \cdot1_{\{Z_{\widetilde{Q}}
 = \int_{\mathfrak{Z}_{\mathcal P} }Z_{P
} \mathrm{d} \widetilde{\lambda} \}}(\omega),
\end{equation}
reminiscent of \textup{(\ref{01})}, where the random variable $  \delta\dvtx
\Omega\rightarrow[0,1] $ is chosen so that \textup{(\ref{testB})} is
satisfied by this $ \widetilde{\varphi} $.
\end{corollary*}

\section{Extensions and ramifications}\label{sec5}

\begin{remark}\label{Remark Assumpions}
The~weak compactness of the set of alternative densities
$\mathfrak{Z}_{\mathcal Q}$ in
Assumption~\ref{Assumption}(i) seems to be crucial. Without it,
we can still get, by Fenchel duality,
\[
\max_{\varphi\in\Phi_{\alpha}} \Bigl(
\inf_{Q\in\mathcal Q}
\mathbb{E}^Q[\varphi]\Bigr)
= \inf_{Q\in\mathcal Q}\Bigl(
\max_{\varphi\in
\Phi_{\alpha}}\mathbb{E}^Q[\varphi]\Bigr),
\]
endowing $\mathbb{L}^\infty$ with the norm
topology. There is no guarantee anymore, however, that the
infimum will be attained in $\mathcal Q$.
The~infimum will be attained at some element $ \widehat{\mu}
  $ of the set $ \mathfrak{M} \subseteq\{\mu\in
\mathfrak{ba}(\Omega,\mathcal{F})_+   \mid  \mu(\Omega)=1\} $,
which contains $  \mathcal{Q}  $; however, a Hahn
decomposition might not exist for this
$\widehat{\mu}$, so we do not generally obtain the 0--1
structure in (\ref{012}) of the
primal solution with respect to the dual solution.

It seems reasonable to endow $\mathbb{L}^\infty$ with the weak-*
topology and to apply Fenchel duality. But it then becomes difficult to
show that a suitable constraint qualification (e.g., that
$ \rho(\varphi)=\sup_{Q\in\mathcal Q}\mathbb{E}^Q[\varphi] $ be
weakly-* continuous at some $\varphi_0\in\Phi_\alpha$) is
satisfied, which is needed to obtain strong duality.

Assuming $\mathfrak{Z}_{\mathcal Q}$ to be weakly compact and
convex, as we have done throughout the present work, has enabled
us to apply a min--max theorem and to ensure that the infimum in
the dual problem is attained within $\mathfrak{Z}_\mathcal
Q\subseteq\mathbb{L}^1$.

\noindent$\bullet$ If we were to drop the weak compactness
Assumption~\ref{Assumption}(ii) on the set of densities
$ \mathfrak{Z}_{\mathcal P} $, then the norm-dual of $\mathcal
L$ would be $ \mathfrak{ba}(\mathfrak{Z}_{\mathcal
P},\mathcal{B}) $ instead of $\Lambda$ (recall the definitions
of the spaces $\mathcal L$ and $\Lambda$ from the start of the
proof of Lemma \ref{theorem opt Test NP}). The~elements of the
space $ \mathfrak{ba}(\mathfrak{Z}_{\mathcal P},\mathcal{B}) $
are the ultimate generalized Bayesian priors: they are allowed to
assign negative weights to sets of possible scenarios and to be
just finitely (as opposed to countably) additive. However, in such a
setting, the Fubini--Tonelli theorem can no longer be applied. It
is possible to endow $\mathcal L$ with the Mackey topology, to
ensure that $\Lambda$ is the topological dual of $ \mathcal L $,
but proving strong duality under this topology is a challenge.
Throughout this paper, Assumption~\ref{Assumption}(ii) is imposed
to ensure that the norm-dual of the space $ \mathcal L $ is
$ \Lambda $ and that a strong duality result can be obtained.

The~weak compactness Assumption 4.1(ii) on
$ \mathfrak{Z}_{\mathcal P}$ can be dropped when using a
different duality approach as in Cvitani\'c and Karatzas
\cite{CviKar01} (cf.~Section~\ref{sec6}), or a utility
maximization duality approach similar to F\"ollmer
and Leukert \cite{FL00}, Section 7. In both cases, however, one would lose
some information about the structure of the optimal randomized
test since $ \widetilde{\varphi} $ will no longer be expressed
with respect to the original family of probability measures
$\mathcal{ P}$, but with respect to some enlarged set (see
Section~\ref{sec6}).
\end{remark}

\begin{remark}
The~results in this paper can be extended in several
directions.
For instance, our proofs have not used the assumption that $ \mathcal
{P} $
and $ \mathcal{Q} $ are families of probability measures. The~results still hold if we instead consider two arbitrary
subsets of $\mathbb{L}^1$, namely, $ \mathcal G $ (in lieu
of $ \mathfrak{Z}_{\mathcal P} $) and~$ \mathcal H $ (in lieu
of $ \mathfrak{Z}_{\mathcal Q} $), that satisfy
Assumption~\ref{Assumption}, as well as
$ \sup_{G\in\mathcal G}\|G\|_{ \mathbb{L}^1}< \infty $.

Furthermore, instead of a constant $\alpha\in(0,1)$, we may
consider a positive continuous functional $ \alpha\dvtx \mathcal
G\rightarrow\mathbb{R}_+  $. The~corresponding optimization problem is
then
%
\begin{equation}\label{testopallg}
\sup_{\varphi\in\Phi} \Bigl( \inf_{H\in\mathcal H}
\mathbb{E}[\varphi H]\Bigr),
\end{equation}
subject to
%
\begin{equation}\label{testopallgNB}
\mathbb{E}^R[\varphi G] \leq
\alpha(G) \qquad   \forall G\in\mathcal G.
\end{equation}
The~problem (\ref{testopallg})--(\ref{testopallgNB}) is no longer
of hypothesis-testing form, in the classical sense, but its
structure is similar to that of testing composite hypotheses. Such
so-called
`generalized hypothesis-testing problems'
also arise in the context of hedging contingent claims in
incomplete or constrained markets, for instance, when one tries to
minimize the expected hedging loss (see Cvitani\'c
\cite{Cvitanic00}, Rudloff \cite{Rudloff07} or, in a related
context, Schied \cite{Schied04}).
This kind of generalized hypothesis-testing problem was studied
for the case of a simple alternative (i.e., $\mathcal H$ being a
singleton) and a positive, bounded and measurable function
$\alpha(\cdot)$, by Witting \cite{Witting85}, Section~2.5.1. For
this case, it was shown with Lagrange duality that the generalized
0--1 structure (\ref{012}) of a test is sufficient for
optimality. Furthermore, it was shown in \cite{Witting85} that
for a finite set $ \mathcal G $, the conditions (\ref{NP Loesung
1}), (\ref{NP Loesung 2}) are necessary and sufficient for
optimality. The~proof of Lemma~\ref{theorem opt Test NP} shows
that a
generalization of these results is even possible
when both the `hypothesis' set $\mathcal G$ and the `alternative
hypothesis' set $\mathcal H$ are infinite, if they satisfy
the above conditions (Assumption~\ref{Assumption} and
$ \sup_{G\in\mathcal G}\|G\|_{ \mathbb{L}^1}< \infty $) and
$ \alpha(\cdot) $ is a positive, continuous function.
\end{remark}

\section{Comparisons and conclusion}\label{sec6}

The~problem of Neyman--Pearson-type testing of a composite null
hypothesis against a
simple alternative has a long history--it has been
considered in a variety of papers and in
several books: Ferguson \cite{Ferguson67}, Witting
\cite{Witting85}, Strasser \cite{Strasser85}, Vajda
\cite{Vajda89}, Lehmann \cite{LehmannRomano05}.

Results on testing a composite hypothesis
against a composite alternative have been obtained in Lehmann \cite
{Lehmann52}, in Krafft
and Witting \cite{KafftWitting67} (which is apparently the first work
to introduce ideas of convex duality to the theory of hypothesis
testing) and in Baumann \cite{Baumann68}, Huber and Strassen~\cite{HubStr73}, \"Osterreicher \cite{Oesterreicher78}, Vajda
\cite{Vajda89}, pages 361--362 and Schied \cite{Schied04}. Lehmann
\cite{Lehmann52} works with a finite set $\mathfrak{Z}_\mathcal{Q}
$, provides existence results and shows
that the composite testing problem can be reduced to one with simple
hypotheses (consisting of the optimal mixed strategy). Krafft and
Witting \cite{KafftWitting67} and Witting \cite{Witting85} use
Lagrange duality and show that even without any compactness
assumptions on the sets $  \mathfrak{Z}_\mathcal{P} $ and $
\mathfrak{Z}_\mathcal{Q} $, the generalized 0--1 structure
(\ref{012}) of $ \widetilde{\varphi} $ is sufficient for
optimality, as well as necessary and sufficient if a dual solution
exists (e.g., when $  \mathfrak{Z}_\mathcal{P} $ and $
\mathfrak{Z}_\mathcal{Q} $ are both finite). In this paper, we show
that under Assumption~\ref{Assumption}, the generalized \mbox{0--1}
structure of (\ref{012}) is necessary and sufficient for optimality,
due to strong duality with respect to the Fenchel dual problem; the
existence of a dual solution then follows from strong duality.
Baumann \cite{Baumann68} proves the existence of a
max--min optimal test using duality results from linear programming and
weak compactness
arguments. The~problem is also studied for densities that are
contents and not necessarily measures.

Huber and Strassen
\cite{HubStr73} dispense with the assumption that all measures in
$  \mathfrak{Z}_\mathcal{P} $ and $
\mathfrak{Z}_\mathcal{Q} $ be absolutely continuous with respect
to a reference measure $R $, at the expense of assuming that
these two sets can be described in terms of ``alternating
capacities'', in the sense of Choquet; for related results, see Rieder
\cite{Rieder77} and Bendarski \cite{Bendarski81,Bendarski82}.

Finally, using totally
different methods and motivated by optimal investment problems in
mathematical finance, Schied \cite{Schied04} studies variational
problems of Neyman--Pearson type for convex risk measures and for
law-invariant robust utility functionals, obtaining explicit
solutions for quantile-based coherent risk measures that satisfy
the Huber--Strassen--Choquet alternating capacity conditions.

One of the most recent works on this subject is the paper by Cvitani\'c
and Karatzas \cite{CviKar01}.
These authors introduce the enlargement
%
\begin{equation}
\label{enlarge}
\mathcal W:= \{W\in\mathbb{L}^1_+ \mid
\mathbb{E}^R[\varphi
W]\leq\alpha, \forall \varphi\in\Phi_\alpha  \}
\supseteq
\co(\mathfrak{Z}_{\mathcal P})
\end{equation}
of the convex hull of the Radon--Nikodym densities of
$ \mathcal{P }$. This `enlarged' set $ \mathcal W $ is convex,
bounded in $\mathbb{L}^1$ and closed under $ R$-a.e.~convergence.
Furthermore, they assume
that the set of densities of $\mathcal Q$ is
convex and closed under $ R$-a.e. convergence. The~starting point
of
\cite{CviKar01} is the observation that
%
\begin{equation}
\label{c/k weak duality}
\forall Q\in\mathcal{Q},\forall  W\in\mathcal W,
\forall z>0, \forall \varphi\in\Phi_\alpha,\qquad
\mathbb{E}^Q[\varphi]\leq\mathbb{E}^R [(Z_Q-zW)^{+}]
+\alpha
z.
\end{equation}
The~existence of a quadruple $(\widehat
Q,\widehat{W},\widehat{z},\widehat\varphi)\in \mathcal
Q\times\mathcal W\times(0,\infty)\times\Phi_\alpha $ which
satisfies (\ref{c/k weak duality}) as equality is then shown and
the structure
%
\begin{equation} \label{opt phi c/k}
\widehat{\varphi}(\omega) = 1_{\{ \widehat{z}  \widehat{W}
< Z_{\widehat{Q}} \} } (\omega) + \delta(\omega) \cdot
1_{\{ \widehat{z}  \widehat{W}
= Z_{\widehat{Q}} \} } (\omega)
\end{equation}
for the optimal randomized test $ \widehat\varphi $ is deduced.
Here,
the triple $ (\widehat
Q,\widehat{W},\widehat{z}) $ is a solution of the optimization
problem
%
\begin{equation}
\label{opt pair c/k}
\mathop{\inf_{z>0}}_{(Q,W) \in  \mathcal Q \times\mathcal W }
\bigl(
\alpha z+\mathbb{E}^R[(Z_{Q}-zW)^+]\bigr)
\end{equation}
and the random variable $ \delta: \Omega\rightarrow[0,1] $ is
chosen so that $  \sup_{P \in\mathcal{P}} \mathbb{E}^{P} [
\widetilde{\varphi}  ] = \alpha $.

\textit{The~methodology of the present paper obviates the need to
introduce the enlargement set $ \mathcal W$ of} (\ref{enlarge}).
Thus, we provide a result about the structure of the solution
$ \widetilde{\varphi} $ in terms of the original families of
probability measures $ \mathcal P  $ and $ \mathcal Q $. We do, however,
need the set $ \mathfrak{Z}_\mathcal Q $ to be weakly
compact.

\noindent$\bullet $ Let us study the relationship between
Theorem~\ref{Theorem Testtheorie} and the results of Cvitani\'c and
Karatzas~\cite{CviKar01}. From the Fubini--Tonelli theorem, it is easy
to
show that
%
\begin{equation}
\label{KL}
k (\lambda)  \int_{\mathfrak{Z}_{\mathcal P}} Z_{P}\, \mathrm{d}\lambda \in
 \mathcal
W \qquad   \forall \lambda\in\Lambda_{+},
\end{equation}
where $k (\lambda) = (\lambda(\mathfrak{Z}_{\mathcal P}))^{-1}$ if
$\lambda(\mathfrak{Z}_{\mathcal P}) > 0$, and $ k (\lambda)=0 $ if
$\lambda( {\mathfrak{Z}_{\mathcal P}})=0$. The~case
$\lambda(\mathfrak{Z}_{\mathcal P})=0$ implies that $\lambda(B)=0$ for
all $B\in\mathcal B$ and thus $\int_{\mathfrak{Z}_{\mathcal P}}
Z_{P }\,\mathrm{d}\lambda=0$. If, in (\ref{c/k weak duality}), we consider only
elements $W$ of the form $ k\int_{\mathfrak{Z}_{\mathcal P}} Z_{P
}\,\mathrm{d}\lambda\in\mathcal W $, then the inequality (\ref{c/k weak
duality}) coincides with the weak
duality between\vspace*{1pt} the primal and dual objective functions
$\mathfrak{p}(Q)$ and $\mathfrak{d}(Q)$, and reduces to the
inequality in (\ref{dual}), whereas
Problem (\ref{opt pair c/k}) reduces to (\ref{opt pair test}).

To summarize, Cvitani\'c and Karatzas
\cite{CviKar01} proved the existence of primal and dual
solutions that satisfy (\ref{c/k weak duality}) as equality. In
order to do this, strong closure
assumptions
had to be imposed. In the methodology of the present paper, the
validity of strong duality, hence also the equality in
(\ref{dual}), are shown directly via Fenchel duality; the
existence of a dual solution then follows. Both methods lead to a
result about the structure of an optimal
test. However, it is now possible, as in Theorem~\ref{Theorem
Testtheorie}, to show the impact of the original
family $ \mathcal P  $ on the sets that define the solution
$ \widehat{\varphi} $ in \cite{CviKar01}, as in (\ref{opt phi
c/k}):
%
\begin{equation}
\widehat{z} \widehat{W} = \int_{\mathfrak{Z}_{\mathcal P}}
Z_{P } \,\mathrm{d}\widetilde{\lambda},
\end{equation}
where $(\widetilde{Q},\widetilde{\lambda})$ is the
pair that attains the infimum in (\ref{opt pair test}), $  \widehat
{z}=\widetilde{\lambda}(
{\mathfrak{Z}_{\mathcal P}}) $ and
\[
\widehat{W} = k (\widetilde{ \lambda})  \int_{\mathfrak
{Z}_{\mathcal P} }Z_{P
} \,\mathrm{d}\widetilde{\lambda},\qquad  \mbox{in the notation of (\ref{KL})}.
\]
No embedding of
$ \mathfrak{Z}_{\mathcal P}  $ into the larger set $\mathcal W$ of
(\ref{enlarge}) is any longer necessary.
However, instead of assuming that
$ \mathfrak{Z}_\mathcal{Q}=\{Z_Q \mid  Q\in\mathcal Q\} $ is
closed under $ R$-a.e. convergence, we need to assume here that
this set is weakly compact in $ \mathbb{L}^1$.

\section*{Acknowledgements}

This research was supported in part by the National Science Foundation,
under Grants DMS-06-01774 and DMS-09-05754.

\printhistory

\end{document}